\documentclass[12pt]{article}
\usepackage{theos}

\def\fileversion{v1.0}
\def\filedate{1995/02/09}
\RequirePackage{amsfonts}

\def\qx{{\ensuremath{\mathbb Q}}}

\newcommand{\p}{\partial}

\newcommand{\ab}{\overline{a}}
\newcommand{\bb}{\overline{b}}
\newcommand{\abx}{\overline{a}(x)}
\newcommand{\bbx}{\overline{b}(x)}
\newcommand{\Rx}{R[x]}
\newcommand{\Rxy}{R[x,y]}
\newcommand{\Rbx}{\overline{R}[x]}
\newcommand{\autRbx}{Aut_{\overline{R}}\overline{R}[x]}
\newcommand{\autRxy}{Aut_RR[x,y]}
\newcommand{\atx}{\tilde{a}(x)}
\newcommand{\btx}{\tilde{b}(x)}

\newcommand{\btf}{\tilde{b}(f_1)}
\newcommand{\kb}{\overline{k}}
\newcommand{\pin}{\frac{1}{p}}

\begin{document}
\title{A note on $k[z]$-automorphisms in two variables}
\author{Eric Edo, Arno van den Essen, Stefan Maubach\footnote{Funded by Veni-grant of council for the
physical sciences, Netherlands Organisation for scientific research (NWO).
Partially funded by the Mathematisches
Forschungsinstitut Oberwolfach as an Oberwolfach-Leibniz-Fellow.}}
\maketitle

\abstract{We prove that for a polynomial $f\in k[x,y,z]$ equivalent are:
(1)$f$ is a $k[z]$-coordinate of $k[z][x,y]$, and (2) $k[x,y,z]/(f)\cong k^{[2]}$ and $f(x,y,a)$ is a coordinate in $k[x,y]$ for some $a\in k$.
This solves a special case of the Abhyankar-Sathaye conjecture.
As a consequence we see that a coordinate $f\in k[x,y,z]$ which is also a $k(z)$-coordinate, is a $k[z]$-coordinate.
We discuss a method for constructing automorphisms of $k[x,y,z]$, and observe that the Nagata automorphism occurs naturally as the first non-trivial
automorphism obtained by this method - essentially linking Nagata with a non-tame $R$-automorphism of $R[x]$, where $R=k[z]/(z^2)$.}

\section*{Introduction}

Let $k$ be a field of characteristic zero. The most famous
$k[z]$-automorphism of $k[x,y,z]$ is undoubtedly Nagata's
automorphism $\sigma:k^3\rightarrow k^3$ given by
$$\sigma=(x-2sy-s^2z,y+sz,z), \mbox{ where } s=xz+y^2$$
In a landmark paper [6] Shestakov and Umirbaev solved the
long standing Nagata Conjecture, asserting that $\sigma$ is
not tame. In January 2007 the authors were rewarded with the
Moore prize for the best research paper in the last six years.
Nagata's automorphism can be constructed in several ways: for example
it was constructed by Nagata in [5] as the composition
${\sigma_1}^{-1}\sigma_2\sigma_1$, where
$\sigma_1=(x+z^{-1}y^2,y)$ and $\sigma_2=(x,y+z^2x)$.

Another construction uses locally nilpotent derivations, namely
one easily verifies that
$\sigma=exp(sD)$, where $D$ is the locally nilpotent derivation
given by $D=-2y\p_x+z\p_y$. Yet another construction can be found
in [2] and [1]. It goes roughly as follows (for a detailed
description we refer to section 1): start with an element
$p$ in $k[z]$, which is no constant, and form the quotient ring
$R=k[z]/(p)$. Consider the polynomial ring in one variable $x$
over $R$ and let $a(x)\in Aut_R R[x]$. To this $R$-automorphism
one constructs a $k[z]$-automorphism $(f_1,f_2)$ of $k[z][x,y]$, which
in turn gives a $k$-automorphism of $k[x,y,z]$. Now the
Nagata automorphism can be found by taking the simplest
non-trivial case in the above construction, namely
 $p=z^2$ and $a(x)=x+\overline{z}x^2$, (see Remark 2.3).

As a consequence of the main result of this paper, Proposition 2.1,
one obtains that in general the $k$-automorphism $(f_1,f_2,z)$ is
tame if and only if the $k[z]$-automorphism $(f_1,f_2)$ is tame if
and only if $a(x)$ is tame in $Aut_RR[x]$, which just means that
$a(x)$ has degree one in $x$. Consequently Nagata's example is non-tame.
The proof of this result obviously uses one of the main results
of [6] which asserts that a $k[z]$-automorphism of $k[x,y,z]$
is tame if and only if it is tame as a $k$-automorphism.

In the last section we give a result on
$k[z]$-coordinates which to our knowledge is new. It asserts
that a polynomial in $k[x,y,z]$ is a $k[z]$-coordinate if
and only if $k[x,y,z]/(f)$ is $k$-isomorphic to a polynomial
ring in two variables over $k$ and $f(x,y,a)$ is a coordinate
in $k[x,y]$ for some $a$ in $k$. Hence this result
proves a special case of the Abhyankar-Sathaye Conjecture and
furthermore it shows that if $f$ is a coordinate in $k[x,y,z]$
which is also a $k(z)$-coordinate, then it is a $k[z]$-coordinate.

\section{Constructing $R$-automorphisms of $R[x,y]$}

In this section we recall a construction of $R$-coordinates ($R$-automorphisms)
which already can be found in [2] and [1].

Let $R$ be a commutative ring and let $p\in R$ be neither a unit nor a
zero-divisor in $R$. Put $\overline{R}=R/Rp$. Let $a(x),b(x)$ in $R[x]$ be such that
$$\ab(\bbx)=x=\bb(\abx) \mbox{ in } \Rbx.$$
\noindent Equivalently
$$a(b(x))=x=b(a(x)) \mbox({mod}p) \mbox{ in }\Rx\,\,\,\,\,\, (1)$$
To such an element $\abx\in \autRbx$ one constructs an element of $\autRxy$
as follows. Put
$$f_1=py+a(x)\,\,\,\,\,\,\, \,\,\,\,\,(2)$$
Since $f_1=a(x)(\mbox{mod}p)$ it follows from (1) that $b(f_1)-x=0$(mod$p$).
Consequently, since $p$ is a non-zero divisor in $R$, there exists
a unique element $f_2$ in $R[x,y]$ such that
$$b(f_1)-x=pf_2\,\,\,\,\,\,\,\,(3)$$

\begin{lem} $F=(f_1,f_2)\in\autRxy$.
\end{lem}

\noindent\textbf{Proof.} Let $g_1=b(x)-py$. Then $a(g_1)=a(b(x))=x$(mod)$p$. So
$$x-a(g_1)=pg_2\,\,\,\,\,\,\,\,(4)$$
for some $g_2\in\Rxy$.

Now we will show that $G=(g_1,g_2)$ is the inverse of $F$: namely by (3)
$$g_1(f_1,f_2)=b(f_1)-pf_2=x\,\,\,\,\,\,\, (5)$$
Furthermore, using (4), (5) and (2) we obtain
$$pg_2(f_1,f_2)=f_1-a(g_1(f_1,f_2))=f_1-a(x)=py.$$
So, since $p$ is no zero divisor in $R$ we get $g_2(f_1,f_2)=y$.

\section{Tame $R$-automorphisms of $\Rxy$}

Let again $R$ be a commutative ring and $n$ a positive integer. An
$R$-automorphism of $R^{[n]}$ is called {\em tame} if it is a finite product
of automorphisms of the form
$$(x_1,\ldots,x_{i-1},ux_i+v(x),x_{i+1},\ldots,x_n)$$
where $u\in R^*$ and $v(x)\in R^{[n]}$ does not contain $x_i$.
The group of tame automorphisms of $R^{[n]}$ is denoted by $T(n,R)$.
So $T(1,R)$ consists of the elements $ux_1+v$ with $u\in R^*$ and $v\in R$
arbitrary.

\medskip

\noindent From now on we assume that $R$ is a {\em domain}.

\medskip

\noindent Keeping the notations from the previous section, the main result
of this section, Proposition 2.1, asserts that if one starts with an element
$\abx\in \autRbx$ and constructs the corresponding element
$F=(f_1,f_2)$ in $\autRxy$, then $F$ is tame if and only if $a(x)$ is tame, in
other words if and only if $a(x)$ is of degree one in $x$. More precisely

\begin{prop} Let $\abx\in \autRbx$ with inverse $\bbx$, $f_1=py+a(x)$ and
$f_2=(b(f_1)-x)/p$. Then there is equivalence between\\
\noindent i) $\abx\in T(1,\overline{R})$.\\
\noindent ii) $(f_1,f_2)\in T(2, R)$.\\
\noindent Furthermore, in this situation $a(x)=a_0+a_1x+px^2\tilde{a}(x)$,
for some $a_0,a_1$ in $R$ and $\atx\in R[x]$ and there exist $c,d$ in $R$
such that $da_1-cp=1$ and $f_2=cx+d(y+x^2\atx)+\btf$ for some $\btx\in R[x]$.
\end{prop}

\noindent\textbf{Proof.} Write $a(x)=\sum_{i=0}^d a_i(x)x^i$ and
$b(x)=\sum_{i=0}^eb_ix^i$. Assume i). Then $\overline{a_1}\in(\overline{R})^*$ and
$\overline{a_i}=0$ for all $i\geq 2$. Since $\overline{a_1}\in(\overline{R})^*$
there exist $c,d\in R$ such that $da_1-cp=1$. Since $\overline{a_i}=0$
for $i\geq 2$ it follows that $p$ divides each such $a_i$ in $R$, so
$$a(x)=a_0+a_1x+px^2\atx \mbox{ for some } \atx\in R[x].$$
Since $da_1=1$(mod $p$) it follows that the inverse of
$\overline{a_0}+\overline{a_1}x$ is equal to $\overline{d}(x-\overline{a_0})$,
whence $\bbx=\overline{d}(x-\overline{a_0})$, so $b(x)=d(x-a_0)+p\btx$,
for some $\btx\in R[x]$. Consequently
$$f_2=(\btf-x)/p=(d(py+a_1x+px^2\atx)+p\btf -x)/p.$$
Using $(da_1-1)x=cpx$  it follows that
$$f_2=cx+d(y+x^2\atx)+\btf.$$
So
$$(f_1,f_2)=(a_1x+p(y+x^2\atx)+a_0,cx+d(y+x^2\atx)+\btf).$$
Now one easily verifies that $(f_1,f_2)\in T(2,R)$, since
$$(x-a_0,y)\circ(x,y-\btx)\circ(f_1,f_2)\circ(x,y-x^2\atx)=(a_1x+py,cx+dy)$$
and $ad-pc=1$. So i) implies ii).\\
Conversely, assume ii). If $\abx\notin T(1,\overline{R})$, then
$d_1:=$deg$_x\abx\geq 2$ and we can write
$$a(x)=\sum_{i=0}^{d_1} a_ix^i+px^{{d_1}+1}\atx$$
for some $\atx$ in $R[x]$ and $p$ does not divide $a_{d_1}$.
Similarly, let $e_1:=$deg$_x\bbx$. Since $\bbx$ is the inverse
of $\abx$ and $d_1\geq 2$, it follows that $e_1\geq 2$. So we can write
$b(x)=\sum_{i=1}^{e_1} b_ix^i+p\btx$, for some $b_i$ in $R$ and
$\btx$ in $R[x]$, where
$p$ does not divide $b_{e_1}$. Then
$$f_1=p(y+x^{{d_1}+1}\atx)+\sum_{i=0}^{d_1} a_ix^i$$
and
$$f_2=\pin[\sum_{j=0}^{e_1}b_j(p(y+x^{{d_1}+1}\atx)+\sum_{i=0}^{d_1}a_ix^i)^j-x]+
f^{{e_1}+1}\btf.$$
Since $(f_1,f_2)\in T(2,R)$ it follows that $(f_1,f_2-f_1^{{e_1}+1}\btf)\in
T(2,R)$ and hence, replacing $y$ by $y-x^{{d_1}+1}\atx$, that
$$(py+\sum_{i=0}^{d_1}a_ix^i,\pin[\sum_{j=0}^{e_1}b_j(py+\sum_{i=0}^{d_1}
a_i x^i)^j-x]\in T(2,R).$$
Since $d_1\geq 2$ the highest degree $xy$-term of $f_1$(resp.$f_2$) equals
$a_{d_1}x^{d_1}$(resp.\\
\noindent $\pin b_{e_1}(a_{d_1}x^{d_1})^{e_1}$).
So by Corollary 5.1.6 of [3], using that $e_1\geq 2$, it follows that there
exists $c$ in $R$ such that
$$\pin b_{e_1}a_{d_1}^{e_1}x^{d_1e_1}=ca_{d_1}^{e_1}x^{d_1e_1}$$
whence $b_{e_1}=pc$, so $p$ divides $b_{e_1}$, a contradiction.
So $\abx\in T(1,\overline{R})$ as desired.

\begin{cor}Let $k$ be a field of characteristic zero and $p\in k[z]$,
but not in $k$. Put $R=k[z]/(p)$.
Let $a(z,x)$ and $b(z,x)$ in $k[z,x]$ be such that
$a(\overline{z},x)\in Aut_RR[x]$ with inverse $b(\overline{z},x)$.
Put $f_1=py+a(z,x)$ and $f_2=(b(z,f_1)-x)/p$. Then $(f_1,f_2,z)\in T(3,k)$
if and only if $(f_1,f_2)\in T(2,k[z])$ if and only if
$a(\overline{z},x)\in T(1,R)$ if and only if
deg$_x a(\overline{z},x)=1$.
\end{cor}

\noindent\textbf{Proof.} Proposition 2.1 gives the equivalence
of the first two statements and
a result of [6] gives the equivalence of the second
and third statement. The last equivalence is obvious.

\medskip

\noindent\textbf{Remark 2.3} If $p=z^2$ and $a(z,x)=x+zx^2$
one obtains the simplest non-tame automorphism, namely
$x+\overline{z}x^2$. The corresponding (non-tame)
$k$-automorphism $(f_1,f_2,z)$ is, apart from a
permutation of the variables $x$ and $y$, the Nagata
automorphism $\sigma$. More precisely
$$\sigma=(y,x,z)\circ (f_1,f_2,z)\circ (y,x,z).$$
So Nagata's automorphism is the simplest non-tame
automorphism ``coming from'' a one dimensional
example.

\section{A remark on $k[z]$-coordinates}

Throughout this section $k$ is a field of characteristic
zero. If $n\geq 1$ and $f$ is a $k$-coordinate in the polynomial ring
$k^{[n]}$ then $k^{[n]}/(f)$ is $k$-isomorphic to a
polynomial ring in $n-1$ variables over $k$. The Abhyankar-Sathaye
Conjecture asserts that the converse is true. This conjecture
is still open for all $n\geq 3$.

In this section we prove a special case of this conjecture in
case $n=3$. More precisely we show

\begin{prop} Let $f\in k[x,y,z]$ be such that $A=k[x,y,z]/(f)$
is $k$-isomorphic to $k^{[2]}$ and that for some $a\in k$
the polynomial $f(x,y,a)$ is a coordinate in $k[x,y]$.
Then $f$ is a $k[z]$-coordinate in $k[x,y,z]$.
\end{prop}

To prove this proposition we need the following result from [4]:

\begin{thm} Let $R$ be a $\qx$-algebra and $f\in R[x,y]$. Let
$D$ be the derivation $f_y\p_x-f_x\p_y$ on $R[x,y]$. Then there is
equivalence between\\
\noindent i) $f$ is a coordinate in $R[x,y]$.\\
\noindent ii) $D$ is locally nilpotent and
$1\in R[x,y]f_x+R[x,y]f_y$.
\end{thm}

\noindent\textbf{Proof of Proposition 3.1} i) Let $\kb$ be an
algebraic closure of $k$ and view $f$ in $\kb[x,y,z]$.
The hypothesis implies that $\kb[x,y,z]/(f)\simeq \kb^{[2]}$
and that $f(x,y,a)$ is a coordinate in $\kb[x,y]$. We will
deduce in ii) below that $f$ is a $\kb[z]$-coordinate
in $\kb[x,y,z]$. It then follows from Theorem 3.2
that $D$ is locally nilpotent on $\kb[x,y,z]$, and hence
on $k[x,y,z]$. Also we obtain from Theorem 3.2 that
 $1\in\kb[x,y,z]f_x+\kb[x,y,z]f_y$,
which implies that $1\in k[x,y,z]f_x+k[x,y,z]f_y$, since
$f$ has coefficients in $k$. Then again applying Theorem 3.2
gives that $f$ is a $k[z]$-coordinate of $k[x,y,z]$. So
we may assume that $k=\kb$.\\
\noindent ii) Now assume $k=\kb$. The hypothesis implies that
$$A/(z-a)\simeq k[x,y,z]/(f,z-a)\simeq k[x,y]/(f(x,y,a))\simeq k^{[1]}.$$
\noindent Since $A\simeq k^{[2]}$ the Abhyankar-Moh theorem implies that
$z-a$ is a coordinate in $A$ and hence so is $z-b$ for all $b\in k$. So
$k[x,y]/(f(x,y,b))\simeq A/(z-b)\simeq k^{[2]}$. Hence by the
Abhyankar-Moh theorem $f(x,y,b)$ is a coordinate in $k[x,y]$ for
all $b\in k$. In particular for all $b$ in $k$ the element $1$ is in the
ideal generated by $f_x(x,y,b)$ and $f_y(x,y,b)$ in $k[x,y]$. Consequently
$f_x(x,y,z)$ and $f_y(x,y,z)$ have no common zero in $k^3$. So by
the Nullstellensatz we have

\medskip

\noindent (6) $1$ is in the ideal generated by $f_x$ and $f_y$ in $k[z][x,y]$.

\medskip

\noindent Now let $D=f_y\p_x-f_x\p_y$ and let $d$ be the maximum of
the $x$ and $y$ degrees of $f_x$ and $f_y$. Since for each $b$ in $k$
the polynomial $f(x,y,b)$ is a coordinate in $k[x,y]$, the derivation
$D$ evaluated at $z=b$ is locally nilpotent and hence it follows from
[3], Theorem 1.3.52 that $D^{d+2}(x)(z=b)=0$ and $D^{d+2}(y)(z=b)=0$
for all $b$ in $k$. This implies that $D^{d+2}(x)=D^{d+2}(y)=0$. So

\medskip

\noindent (7) $D$ is locally nilpotent on $k[z][x,y]$.

\medskip

\noindent Then it follows from (6), (7) and Theorem 3.2 that
$f$ is a $k[z]$-coordinate of $k[x,y,z]$, as desired.

\section*{References}

\noindent [1] J. Berson, Stable tame coordinates, J. Pure and Applied
Algebra, 170 (2002), 131-143.

\medskip

\noindent [2] E. Edo and S. V\'en\'ereau, Length 2 variables and transfer,
Annales polonici Math. 76 (2001), 67-76.

\medskip

\noindent [3] A. van den Essen, Polynomial Automorphisms and the Jacobian
Conjecture, Progress in Mathematics, vol. 190, Birkh\"auser-Verlag, (2000).

\medskip

\noindent [4] A. van den Essen and P. van Rossum, Coordinates in
two variables over a $\qx$-algebra, Trans. of the AMS, vol. 356, no. 5,
(2004), 1691-1703.

\medskip

\noindent [5] M. Nagata, On the automorphism group of $k[X,Y]$, Kyoto
Univ. Lectures in Math. 5 (1972).

\medskip

\noindent [6]  I. Shestakov and U. Umirbaev, The tame and wild automorphisms
of polynomial rings in three variables, J. of the AMS, 17 (2004), no. 1,
197-227.\\

\medskip

\noindent Author's addresses:\\

\noindent
Eric Edo\\
edo@univ-nc.nc\\
\noindent
Porte S 20, Nouville Banian,\\
     Université de Nouvelle Calédonie, BP R4\\
     98 851 Nouméa Cedex, Nouvelle Calédonie.\\

\noindent
Arno van den Essen\\
essen@math.ru.nl\\
\noindent
Faculty  of Science, Mathematics and Computer Science,\\
Radboud University Nijmegen\\
Postbus 9010, 6500 GL Nijmegen \\
The Netherlands\\

\noindent
Stefan Maubach\\
s.maubach@science.ru.nl\\
\noindent
Faculty  of Science, Mathematics and Computer Science,\\
Radboud University Nijmegen\\
Postbus 9010, 6500 GL Nijmegen \\
The Netherlands

\end{document}